\documentclass[11pt]{article}
\usepackage{amsfonts,latexsym}

\title{Parallel transport on higher loop spaces}
\author{Ivan Emilov Horozov}
\date{June 20, 2012}

\newcommand{\beq}{\begin{equation}}
\newcommand{\eeq}{\end{equation}}
\newcommand{\beqa}{\begin{eqnarray}}
\newcommand{\eeqa}{\end{eqnarray}}
\newcommand{\beaa}{\begin{eqnarray*}}
\newcommand{\ben}{\begin{eqnarray*}}
\newcommand{\eaa}{\end{eqnarray*}}
\newcommand{\een}{\end{eqnarray*}}

\newcommand{\text}{\textrm}
\newcommand \nc {\newcommand}
\nc \proof {\noindent {\em{Proof.\/ }}} \nc \qed {$\Box$\hfill}
\newtheorem{theorem}{Theorem}[section]
\newtheorem{lemma}[theorem]{Lemma}
\newtheorem{proposition}[theorem]{Proposition}
\newtheorem{corollary}[theorem]{Corollary}
\newtheorem{definition}[theorem]{Definition}
\newtheorem{example}[theorem]{Example}
\newtheorem{remark}[theorem]{Remark}
\newtheorem{conjecture}[theorem]{Conjecture}
\newtheorem{question}[theorem]{Question}
\nc \bth[1] {\begin{theorem}\label{t#1} } \nc \ble[1]
{\begin{lemma}\label{l#1} } \nc \bpr[1]
{\begin{proposition}\label{p#1} } \nc \bco[1]
{\begin{corollary}\label{c#1} } \nc \bde[1]
{\begin{definition}\label{d#1}\rm } \nc \bex[1]
{\begin{example}\label{e#1}\rm } \nc \bre[1]
{\begin{remark}\label{r#1}\rm } \nc \bcon[1]
{\begin{conjecture}\label{con#1}\rm } \nc \bque[1]
{\begin{question}\label{que#1}\rm }
\nc {\eth} { \end{theorem} } \nc {\ele} { \end{lemma} } \nc
{\epr}{ \end{proposition} } \nc {\eco} { \end{corollary} } \nc
{\ede} {\end{definition} } \nc {\eex} { \end{example} } \nc {\ere}
{\end{remark} } \nc {\econ} { \end{conjecture} } \nc {\eque}
{\end{question} }
\nc \eqref[1] {{\rm{(\ref{#1})}}} \nc \thref[1]{Theorem \ref{t#1}}
\nc \leref[1]{Lemma \ref{l#1}} \nc \prref[1]{Proposition
\ref{p#1}} \nc \coref[1]{Corollary \ref{c#1}} \nc
\deref[1]{Definition \ref{d#1}} \nc \exref[1]{Example \ref{e#1}}
\nc \reref[1]{Remark \ref{r#1}} \nc
\conref[1]{Conjecture\ref{con#1}}

\def \R {{\mathcal R}}

\def \Z {{\mathbb Z}}

\def \R {{\mathbb R}}
\def \C {{\mathbb C}}

\nc \Wr {Wr} \nc \GRN { \Gr^{(N)} }
\nc \GRA[1] { \Gr_A^{(#1)} }   
\nc \GRAN { \GRA{N} } \nc \GrA[1] { \Gr_A(#1) }\nc \GrAa {
\GrA{\alpha} }
\nc \GRB[1] { \Gr_B^{(#1)} }   
\nc \GRBN { \GRB{N} } \nc \GrB[1] { \Gr_B(#1) } \nc \GrBb {
\GrB{\beta} }
\nc \GRMB[1] { \Gr_{MB}^{(#1)} }   
\nc \GRMBN { \GRMB{N} } \nc \GrMB[1] { \Gr_{MB}(#1) } \nc \GrMBb {
\GrMB{\beta} }

\begin{document}

\title{{\LARGE\bf{Parallel Transport on Higher Loop Spaces}}}

\author{
Ivan ~Horozov
}
\date{}
\maketitle

\begin{abstract}
We construct a parallel transport on higher loop spaces of a manifold in term of
a higher dimensional generalization of iterated path integrals.
Under mild assumptions, we define a de Rham complex on higher loop
spaces and we recover a known result \cite{Ha1} and \cite{Ha2} of a de Rham
structure on higher homotopy groups of a manifold.

The key ingredient is a new definition of iterated integrals on
membranes, which also have applications in number theory
\cite{MDZF}, algebraic geometry \cite{Parshin} and mathematical
physics \cite{Physics}.
\end{abstract}

\tableofcontents

\newpage
\section{Introduction}

There has been different attempts to construct higher dimensional
analogues of Chen's iterated integrals of differential forms over
paths. We propose a definition, which is suitable for several
different areas: algebraic and differential topology \cite{DH}, algebraic
geometry \cite{Parshin}, number theory \cite{MDZF}, mathematical physics \cite{Physics}. In this paper we
will be mostly interested in the topological aspects of the
integrals. We are going to construct a sort of differential
geometry on higher order path spaces.

The goal of this paper if to define differential forms on higher
loop spaces and a notion of a parallel transport on higher loop
spaces. As a consequence, we construct de Rham complex on higher
loop space, under mild assumptions.

Chen has used iterated path integrals as functions on path spaces.
In this paper we define higher dimensional analogues, which we
call iterated integrals over membranes. They have a natural meaning
of functions on higher path spaces.

In a previous paper (see \cite{DH}), we have explored functions on
higher order loop spaces in terms of iterated integrals over
membranes, especially, homotopy invariant functions on higher
order loop spaces.

Why are we interested in parallel transport on higher order path spaces?

As a starting point, consider a manifold $M$ and with a connection.
We are going to interpret a parallel transport along a path $\gamma$
in the spirit of Chen. Then the parallel transport can be expressed as
a generating series of iterated integrals of the $1$-form defining the connection
over the path $\gamma$. One can think of an iterated integral over a path $\gamma$ as a function on
the path space
$$PM=Maps([0,1],M).$$

Similarly, we can try to do parallel transport of
a path $\gamma_0$ to $\gamma_1$ along a membrane
$$g:[0,1]^2\rightarrow M,$$
such that $\gamma_0$ is one side of $g$ and
$\gamma_1$ is the opposite side, formally,
$$\gamma_0(t)=g(0,t)$$
and
$$\gamma_1(t)=g(1,t).$$
We want to parallel transport $\gamma_0$ along $g$
in order to obtain information at the path $\gamma_1$.
One can think of $g$ as a function on the $2$-path space
$$P^2M=Maps([0,1]^2,M).$$

Using lowest order of iteration on a membrane $g$ of a suitable, we
express the Riemann curvature tensor.
This idea is developed in \cite{Physics},
where we consider a minor modification of the Riemann curvature tensor
in order to include spinors.

A connection on the loop space (in terms of gerbes) has been used
by Brylinski and MacLaughlin to construct the Parshin symbol on a
surface and prove reciprocity laws for the symbol \cite{BrMcL}.
Similarly, in a separate paper, see \cite{Parshin}, we use iterated
integrals over $2$-dimensional membranes for the purpose of the
Parshin symbol. Moreover, we construct a new symbol, whose
symmetries are the same as the symmetries of the Riemann curvature
tensor. This coincidence of symmetries has a reason: Both the
Riemann curvature tensor and the new symbol of $4$ rational
functions on an algebraic surface can be written as an iterated
integral over a $2$-dimensional membrane. And the multiplication
among these integrals has the meaning of parallel transport.

More generally, we define
{\bf{a parallel transport}} on the $(n-1)$-path space,
$$P^{(n-1)}M=Maps([0,1]^{(n-1)},M)$$
as a iterated integral over an $n$-dimensional membrane, which is a function on the $n$-path space, $$P^nM.$$

Besides making a major step towards de Rham structure on higher
loop spaces, our approach gives an alternative approach to gerbes,
which is categorical way for dealing with connections on loop
spaces. The theory of gerbes has started from Grothendieck and
Giraud \cite{Gi}. Higher analogies were developed by Lurie
\cite{L}.

Let us explain the key steps in finding cohomology of higher loop
spaces. Given a topological space (more precisely, a
differentiable space, see Definition 2.3) Chen considered the loop space as the space
of continuous and piecewise smooth maps from a circle $S^1$ with a
base point $*$ to a differentiable space $X$ with a base point
$x_1$. We use the notation
$$L_{x_1}X=Maps((S^1,*),(X,x_1))$$
for the loop space with compact-open topology.
Let $P_{x_1}X$ be the path space
$$P_{x_1}X=Maps([0,1],(X,x_1)),$$
such that $0$ is sent to $x_1$. Let also $UX$ denote the universal
cover of a differentiable space $X$. Using iterated path integrals
of differential forms on $X$, Chen gave a de Rham structure on a
quotient $\pi_1(X)$ and on $L_{x_1}X$, when $X$ is simply
connected, or equivalently to de Rham structure on $L_{x_1}UX$. In
this paper we define iterated integrals on membranes, which can be
interpreted as iterated integrals on differenital forms on
(higher) loop spaces. Then using Chen's approach, under mild assumptions,  we can give de
Rham structure on a quotient of $\pi_1(L_{x_1}UX)=\pi_2(UX)=\pi_2(X)$, using
$2$-dimensional iterated integrals. More generally, using
$n$-dimensional iterated integrals, we give a de Rham  structure
on a quotient of $\pi_1(L_{x_n}U\dots L_{x_1}UX)=\dots=\pi_n(X)$ and de Rham
complex on $L_{x_{n-1}}U\dots L_{x_0}UX$, which has the homotopy
type of $L^nX=Maps((S^n,*),(X,{x_1}))$. There are technical
conditions that the manifold $X$ satisfy in order to make the
above conclusions. For example, the cohomology of the $(n-1)$-th loop
space $H^*(L_{x_{n-1}}U\dots L_{x_1}UX)$ has to be of finite type.

For a more categorical construction consult \cite{FV}.
The authors have considered an a $n$-fold monoidal categories for
the study of the $n$-th loop spaces. In an analology, in this paper we consider
iterated integrals over $n$-dimensional membranes in order to examine $n$-fold loop spaces. 
A geometric construction for the higher homotopy groups is considered in
\cite{Dev}, where the authors  evaluate Hopf invariants on iterated Whitehead products in
terms of the 'conﬁguration pairing.'



{\bf{Structure of the paper:}} In Section 2, recall some of Chen's
work such as iterated path integrals, shuffle relations, parallel
transport via iterated integrals, de Rham structure on $\pi_1(X)$
and de Rham complex on the loop space of a manifold.

In Section 3, we define higher dimensional analogue of iterated path integrals,
 which we call iterated integrals over membranes.
We use them to define a connection on higher order path spaces.

In Subsection 3.1, we give the the definition of iterated integrals over an
$n$-dimensional membrane. It is essentially $n$ iterated integrals in
$n$ independent directions. It has a meaning in special relativity
as having $n$ observers.

In Subsection 3.2, we define differential forms on higher path
spaces and higher loop spaces, using iterated integrals on
membranes, generalizing Chen's construction to higher loop spaces.

In Subsection 3.3, we define define and prove shuffle relations
analogues to the shuffle relations for iterated path integrals.
The shuffle relations in Proposition 3.1 are slightly more general
than the one considered in \cite{DH}. They play a key role in the
parallel transport on higher order path spaces: In direction of
the parallel transport, one considers cutting as in Subsection
2.7. In all other directions one considers shuffle relations.

In Subsection 3.4, we define a product relation of two integrals
over two membranes, having a common face. This section together with
the section on shuffle relations for iterated integrals on
membranes are essential in the next Subsection.

In Subsection 3.5, we formulate the parallel transport on higher
path spaces in terms of iterated integrals on membranes, using
product relation in one direction and shuffle relations in the
remaining directions.

In Subsection 3.6, we provide a de Rham structure on higher homotopy
groups, under mild assumptions.

In Subsection 3.7 we  construct a de Rham complex on higher loop
spaces, based on Subsections 3.2 and 3.5.

In Section 4, we present a Conjecture that iterated integrals on
membrane capture rational homotopy type together with the
Postnikov tower tensored with the rational numbers.

\section{Iterated integrals over paths}

\subsection{Iterated integrals of $1$-forms over paths (overview)}

Let $\omega_1,\dots,\omega_n$ be differential $1$-forms on $X$. Let $\gamma$
$$\gamma:[0,1]\rightarrow X.$$
be a path on $X.$

\begin{definition}
An iterated integral on $1$-forms over a path is defined as
$$\int_\gamma \omega_1\dots\omega_n
=
\int\dots\int_{0<t_1<\dots<t_n<1}
\gamma^*\omega_1(t_1)\wedge\gamma^*\omega_n(t_n).$$

Let $A_1,\dots,A_n$ be formal non-commuting variable, which
commute with differentiation. (One can think of them as any
constant square matrices of the same size.) Consider the
connection
$$\nabla=d-\sum_{i=1}^nA_i\omega_i.$$
A parallel transport along $\gamma$ with respect to the connection $\nabla$ is given by the generating series
$$\Psi_\gamma=1+\sum_i A_i\int_\gamma\omega_i+\sum_{i,j}A_iA_j\int_{\gamma}\omega_i\cdot\omega_j+\dots.$$

\end{definition}

\begin{theorem}
Let $\gamma_1$ and $\gamma_2$ be two paths such that the end of
$\gamma_1$ is the beginning of $\gamma_2$. Then, we have that the
composition of generating series of iterated integrals over two paths is equal to the generation series over the composition of the two paths. That is,
$$\Psi_{\gamma_1}\Psi_{\gamma_2}=\Psi_{\gamma_1\gamma_2}.$$
\end{theorem}

\subsection{Differential forms on a path spaces}

\begin{definition} (Differentiable space, Chen)
A differentiable space $X$ is a Hausdorff space equipped with a
family of maps called plots, which satisfy the following
conditions:

(a) Every plot is a continuous map of the type $\phi:U\rightarrow
X$, where $U$ is a convex region in $\R^n$ for some $n$.

(b) If $U'$ is also a convex region (not necessarily of the same dimension as $U$) and
if $i:U'\rightarrow U$ is a $C^\infty$ map, then $\phi i$ is also a plot.

(c) Each map $\{0\}\rightarrow X$ is a plot.
\end{definition}

We say that a plot $\psi':U'\rightarrow X$ goes through another plot if $\psi:U\rightarrow X$ via $i$
if there exists a $\C^\infty$ map  $i:U'\rightarrow U$ such that $\phi'=\phi i$.

\begin{definition} (forms on a differentiable space, Chen)
A  $p$-form $\theta$ on a differentiable space $X$ is a rule which
assigns to each plot $\phi:U\rightarrow X$ a $p$-forms
$\theta_\phi$ on the convex region $U$ such that, if $\phi'$ is a
plot which goes though $\phi$ via $i$, then
$$\theta_{\phi'}=i^* \theta_\phi.$$
\end{definition}

\begin{definition}
A piecewise smooth path on a differentiable space $X$ is a
continuous map $\alpha:[0,1]\rightarrow X$ such that, for some
partition $0=t_0<\dots<t_r=1$ of the unit interval, each
restriction $\alpha|[t_{i-1},t_{i}]$ is a plot of $X$.
\end{definition}

Let $P(X)$ denote the space of all piecewise smooth paths on $X$
with the compact open topology. Every map $\alpha:U\rightarrow
P(X)$ gives rise to a map
$$\phi_\alpha:U\times I\rightarrow X$$
given by $(\xi,t)\mapsto \alpha(\xi)(t)$. A plot on $P(X)$ is defined to be a continuous map
$\alpha:U\rightarrow P(X)$, $U$ being s convex region, such that, for some partition
 $0=t_0<\dots<t_r=1$
of the unit interval, the restriction of $\phi_\alpha$ to each
$U\times [t_{i-1},t_i]$ is a plot of $X$.

Let $\Lambda^p(U)$ denotes the space of $p$-forms on the space $U$. Denote also by $\Lambda^*(U)=\oplus_p\Lambda^p(U)$. 

A $\Lambda^p(U)$-valued function on $P(X)$ is a $p$-form in $\Lambda^*(U\times [0,1])$ of the type
$$\theta(t)=\sum\alpha_{i_1,\dots,i_p}(\xi,t)d\xi^{i_1}\wedge\dots\wedge d\xi^{i_p}.$$

For $a,b\in[0,1]$ we define
$$\int_a^b\theta(t)dt
=
\sum\left(\int_a^b\alpha_{i_1,\dots,i_p}(\xi,t)dt\right)d\xi^{i_1}\wedge\dots\wedge d\xi^{i_p},$$
which is again a $p$-valued form of $t$.

Every $(p+1)$-form $\theta$ on $U\times[0,1]$ can uniquely be
written as $dt\wedge \theta'(t)+\theta''(t),$ where $\theta'(t)$
and $\theta''(t)$ are respectively $\Lambda^p(U)$-valued and
$\Lambda^{p+1}(U)$-valued functions of $t$ on $[0,1]$.

Let $\theta$ be a $(p+1)$-form on
a differentiable space $X$. If $\alpha:U\rightarrow P(X)$ is a plot, then the
$(p+1)$-form $\theta_{\phi_\alpha},$
which is piecewise defined on $U\times [0,1]$, can be uniquely written as
$dt\wedge \theta'(t)+\theta''(t)$,
where $\theta'(t)$ and $\theta''(t)$
are piecewise defined $\Lambda(U)$ valued functions of
$t$. We define
$$\theta(\alpha,\dot\alpha)=w'(t),$$
which is a piecewise defined $\Lambda^p(U)$-valued function of $t$ on $[0,1]$.

We write $\Lambda^+(X)=\sum_{p>0}\Lambda^p(X).$

If $\theta_1,\dots \theta_r\in \Lambda^+(X)$, define $\int
\theta_1\dots \theta_r\in \Lambda^*(P(X))$ such that, for any plot
$\alpha:U\rightarrow P(X)$,
$$\left( \int \theta_1\dots \theta_r\right)_\alpha
=
\int_0^1\theta_1(\alpha,\dot\alpha)dt\dots \theta_r(\alpha,\dot\alpha)dt\in\Lambda^*(U).$$

Let $$L_xX=Maps((S^1,*),(X,x))$$ be the space of piecewise smooth maps with a base point $x\in X$.
Let
$$LX=Maps((S^1,*),X),$$ be the total space of loops, whose fiber over a point $x\in X$ is $L_xX$.

\begin{lemma}
In $\Lambda^*(L_x X)$, we have
$$d\int \theta=-\int d \theta$$
and for $r>1$
\begin{eqnarray}
d\int \theta_1\dots \theta_r
=\sum_{i=1}^r (-1)^i \int J\theta_1\dots J\theta_{i-1}(d\theta_i)\theta_{i+1}\dots \theta_r-\\
-\sum_{i=1}^{r-1}(-1)^i\int J\theta_1 \dots (J\theta_1\wedge \theta_{i+1})\dots \theta_r.
\end{eqnarray}
\end{lemma}

\subsection{Shuffle relations}

Let $\theta_1,\theta_2,\dots,\theta_{m'}\dots,\theta_{m'+m''}$
be differential forms on $X$. Let
$\gamma$ be a path on $X$.
\begin{definition}
\label{def sh}
Let $Sh(m',m'')$ be the set of permutations $\rho$ of the set $\{1,2\dots,m'+m''\}$
 such that
$$\rho(1)<\dots<\rho(m')\mbox{ and }\rho(m'+1)<\dots<\rho(m'+m'').$$
\end{definition}

\begin{theorem}
\label{thm shuffle}
We have
$$\int_{\gamma}\theta_1\cdots\theta_{m'}
\int_{\gamma}\theta_{m'+1}\cdots\theta_{m'+m''}
=
\sum_{\rho\in Sh(m',m'')}\int_{\gamma}\theta_{\rho(1)}\cdots\theta_{\rho(m'+m'')}.$$
\end{theorem}

\begin{definition}
\label{def shuffle 1}
Let $\overline{Sh}(m',m'')$ be the set of permutations $\rho$ of the set \\
$\{0,1,\dots,m'+m''+1\}$
 such that $\rho(0)=0$, $\rho(m'+m''+1)=m'+m''+1$ and
$$\rho(1)<\dots<\rho(m')\mbox{ and }\rho(m'+1)<\dots<\rho(m'+m'').$$
\end{definition}

\begin{theorem}
\label{thm shuffle 1}
We have
\begin{eqnarray}
\left[\theta'_{0}(\gamma(0))\left(\int_{\gamma}\theta_1\cdots\theta_{m'}\right)\theta'_{m'+1}(\gamma(1))\right]
\cdot\\
\cdot
\left[\theta''_{0}(\gamma(0))
\left(\int_{\gamma}\theta_{m'+1}\cdots\theta_{m'+m''}\right)\theta''_{m'+m''+1}(\gamma(1))\right]
=\\
=(\theta'_0\wedge\theta''_0)(\gamma(0))
\left(\sum_{\rho\in Sh(m',m'')}\int_{\gamma}\theta_{\rho(1)}\cdots\theta_{\rho(m'+m'')}\right)
(\theta'_{m'+1}\wedge\theta''_{m'+m''+1})(\gamma(1)).
\end{eqnarray}
\end{theorem}

\subsection{Parallel transport along a path}
The purpose of this Subsection is to make a relation between
iterated integrals and parallel transport. The well accepted
notion of parallel transport can be expresses as a generating
series of iterated path integrals of certain type. In some sense,
which we will explain, each iterated integral can be considered as
portion of parallel transport. It will be useful to consider such
a portion of the parallel transport in order to make the analogous
construction for higher loop spaces in the next Section.

First, we consider a parallel transport along a path $\gamma$,
using a connection $\nabla=d-\theta$ of a form $w$.
It is given by
$$w\mapsto \left(\sum_{n=0}^\infty\int_\gamma w\circ\theta^{\circ n}\right)\theta(\gamma(1)).$$

At a finite step the parallel transport is just one of the iterated integrals under the sum. Namely,

$$w\mapsto \left(\int_\gamma w\circ\theta^{\circ n}\right)\theta(\gamma(1)).$$

\subsection{De Rham structure on $\pi_1(X)$}

Let $\theta_i$, $\theta_{ij}$, $\theta_{ijk}$, $\dots$ be
differential forms on $X$. Let $X_i$, $X_{ij}$, $X_{ijk}$, $\dots$
be formal indeterminates. Let
$$\omega=\omega_1+\omega_2+\dots$$
be a formal connection, where
$\omega_1=\sum\theta_iX_i,$ $\omega_2=\sum \theta_i\theta_{j}X_iX_j+\theta_{ij}X_{ij}$, $\dots$. Let
$$k=k_1+k_2+\dots,$$
be the curvature
where
$$k_r=d\omega_r-\sum_{i=1}^{r-1}J\omega_i\wedge \omega_{r-i}.$$

Let $L(X)$ denote the Lie algebra generated formally by
$X_1,\dots,X_m$.

One of the main results of Chen is the following.

\begin{theorem}(Chen, \cite{Ch})
\label{thm pi1} Let $X$ be a differentiable space. Let $W$ be a
vector space of closed $1$-forms on $X$, representing cohomology
classes, with basis $\{w_1,\dots,w_m\}$. Let $V=\{v_1,\dots,v_l\}$
be a basis for $W\wedge W\subset \Lambda^2(X)$. Let the coefficients $c_{ijk}$ be such that
$$w_i\wedge w_j=\sum c_{ijk}v_k.$$
Then there exists an epimorphism of Lie algebras
$$gr \pi_1(X)\otimes_\Z\R\rightarrow L(X)/(N),$$
where $(N)$ is the ideal generated by $\sum_{i,j}c_{ijk}[X_i,X_j]$ for $k=1,\dots,l$.
\end{theorem}

Chen has obtained the above theorem by examining iterated path integrals (see \cite{Ch}). In \cite{Ch2}, he describes the kernel of this map.

\subsection{De Rham complex on loop spaces}

Let $X$ be a path connected differentiable
space. Let $A$ be a differential graded subalgebra of
$\Lambda^*(X)$ such that $dA^0=A^1\cap \Lambda^0(X).$ Let $A'$ be
the subcomplex of $\Lambda^*(L_x X)$ spanned by iterated integrals
of the type $\int \theta_1 \dots \theta_r\in A^+.$ Let
$\hat{C}_*(X)$ be the chain complex of those smooth singular
simplices of $X$ that map the $1$-skeleton of the simplex to the
base point $x_0$ of $X$.

\begin{theorem} (Chen)
\label{thm 1loop} 
With the above notation, if the following conditions hold:

(a) as a topological space, $X$ is simply connected, and its singular homology
is of finite type;

(b) the canonical map from $\hat{C}_*(X)$ into the normalized
singular simplicial chain complex of $X$ is a chain equivalence;

(c) $H(A)\equiv H^*(X;k)$ via  $\hat{C}_*(X)$,

then $H(A')\equiv H^*(L_x X;k).$
\end{theorem}

Lift a differential form $\theta$ on $X$ to a differential form
$\tilde{\theta}$ on $LX$ by the pull-back of $LX\rightarrow X$.
Multiplying (as differentially graded algebra) an iterated
integral over a loop by $\tilde{\theta}$ is
$$\tilde{\theta}_x\left(\int\theta_1\dots\theta_r\right)|L_xX,$$
where $\tilde{\theta}_x$
 is the differential form at the point $x$ and $\sigma_x$ is a loop based at the point $x$.
Globally on $LX$
it will be written as
$$w=\theta\int\theta_1\dots\theta_r.$$
We define
$$(d w)_x
=
(d\theta)_x\left(\int\theta_1\dots\theta_r\right)|L_xX
+
(J\theta)_x \left(d\int\theta_1\dots\theta_r\right)|L_xX.$$

Let $X$ be a path connected differentiable space. Let $A$ be a
differential graded subalgebra of $\Lambda^*(X)$ such that
$dA^0=A^1\cap \Lambda^0(X).$ Let $LA$ be the subcomplex of
$\Lambda^*(LX)$ spanned by iterated integrals of the type
$\theta\int \theta_1 \dots \theta_r\in A^+.$
 \begin{corollary} If the following
conditions hold:

(a) as a topological space, $X$ is simply connected, and its singular homology
is of finite type;

(b) the canonical map from $\hat{C}_*(X)$ into the normalized
singular simplicial chain complex of $X$ is a chain equivalence;

(c) $H(A)\equiv H^*(X;k)$ via  $\hat{C}_*(X)$, then $H(LA)\equiv H^*(LX;k).$
\end{corollary}
\proof It is a direct consequence of Chen's Theorem and of the
fibration
$$L_xX\rightarrow LX\rightarrow X.$$

\section{Iterated integrals over membranes}

\subsection{Definition of iterated integrals on membranes}
Consider an $n$-dimensional unit cube. For each of the directions $i=1,2,\dots,n$ of
the edges, cut such an edge by hyperplanes parallel to the one of
the faces of the cube. Let $t_i^1,t_i^2,\dots$ be variables in the
unit interval, corresponding to points inside an edge in
$i$-direction cut by the hyperplanes. Consider all such $t_i^{j_i}$'s for
$i=1,\dots,n$ and $j_i=0,1\dots,k_i.$

Formally, let
$$D=\{t_i^{j_i}|0=t_i^0<t_i^1<t_i^2<\cdots<t_i^{k_i+1}=1;i=1,\dots,n;j_i=0,\dots,k_i+1\}.$$
Let ${\underline{\bf{k}}_i}=\{0,1,\dots,k_i,+1\}$ and ${\bf{n}}=\{1,2\dots,n\}$.
Let
$$j=(j_1,\dots,j_n)$$
be an $n$-tuple of integers such that $0\leq j_i\leq k_i+1$ for $i=1,\dots,n.$
For each such $j$, let $\omega_j$is a differential form.
Let
$\alpha_j=g^*\omega_j$
Let $J_j\subset {\bf{n}}$.
We define $\alpha_j(J_j;t_1,\dots,t_n)$ as
$$\alpha_j=\sum_{J_j\subset {\bf{n}}}\alpha_j(J_j;t_1,\dots,t_n)\bigwedge_{i\in J_j}dt_i.$$

Let
$$\beta_j=\alpha_j\left(J_j;t_1^{j_1},\dots,t_n^{j_n}\right)\bigwedge_{i\in J_j}dt_i^{\rho_i(j)}.$$

We assume the following condition:

(*) Consider the finite union
$$\bigcup_j\bigcup_{i\in J_j}\{t_i^{j_i}\}.$$
We want that the above union is disjoint and moreover,
 for each $I=1,\dots,n$ and $j_i=1,\dots,k_i$, the element $t_i^{j_i}$ occurs exactly once.

When condition (*) is satisfied, we define an iterated integral over a membrane as
$$\int_{g}^{k,J}\omega_0\dots\omega_k=\int_D\bigwedge_j\beta_j.$$

Due to condition (*), we have that each differential form $dt_i^{j_i}$ occurs
{\it{exactly once under the integral}}.

\subsection{Differential forms on higher loop spaces}
\label{subsec forms on higher loops}
We are going to define
differential forms on higher loop spaces of a manifold in two
ways. First, we will give a definition relating differential forms
$\Lambda$ on a manifold $X$ with differential forms $L^n\Lambda$
on $L^nX$. In the case when $n=1$ and $X$ is a differentiable
space, it is done by Chen, (and recalled in Subsection 2.2). Then
we are going to prove that $L^n\Lambda=L^1(L^{n-1}\Lambda)$.

Let $J_j$ be a subset of $\{1,2\dots,n\}.$
Let $\omega$ be a $(p+|J_j|)$-form on $U\times [0,1]^n$.
Then there exists unique forms $\omega'$ and $\omega''$,
respectively in $\Lambda^{p+c_j}(U)$ and $\Lambda^{p+|J_j|}(U\times[0,1]^n)$,
such that
$$\omega=\left(\bigwedge_{i\in J_j-\rho_i^{-1}(0)-\rho_i^{-1}(k_i)}dt_i\right)\wedge \omega'+\omega'',$$
where $\rho_i^{-1}(0)$ and $\rho_i^{-1}(k_i)$ give boundary points
in direction $i$ and $c_j$ is (codimension of a face of the
$n$-cube). Formally, $c_j$ is the number of $i$'s where
$\rho_i^{-1}(0)\in J_j$ or $\rho_i^{-1}(k_i)\in J_j$. Let
$\alpha:U\rightarrow P^nX$ be a plot. Let $\omega_\alpha$ be
$(p+|J|)$-form on $X$. Express it as
$$\omega_\alpha=\left(\bigwedge_{i\in J_j-\rho_i^{-1}(0)-\rho_i^{-1}(k_i)}dt_i\right)\wedge \omega'+\omega''$$
Let
$$\omega(\alpha,\dot\alpha)=\omega'$$
We define a differential form on $P^nX$ associated to a plot $\alpha$ on $P^nX$ to be
$$\int_{g}^{\rho,k,J}\omega_1(\alpha,\dot\alpha)\dots\omega_m(\alpha,\dot\alpha),$$
where the degree of $\omega_j(\alpha,\dot\alpha)$ is at least $|J_j|-c_j$.

\subsection{Shuffle product}

Let $k'=(k'_1,\dots,k'_n)$, $k''=(k''_1,\dots,k''_n)$ and $j=(j_1,\dots,j_n)$.
Let us define
$$
Sh(k',k'')=\prod_{i=1}^n Sh(k'_i,k''_i),
$$
where the product is over all $i=1,\dots,n$, of shuffles
$\overline{Sh}(k'_i,k''_i)$ as in Definition \ref{def sh}

For $j=(j_1,\dots,j_n)$,
let
$$\omega_j=\omega'_j,$$ and
$J_j=J'_j$
for $1\leq j_i\leq k'_i$ for all $i=1,\dots,n$. 
And let
$$\omega_j=\omega''_{j-k'}$$
and $J_j=J''_{j-k''},$ when
$k'_i<j_i\leq k'_i+k''_i$ for all $i=1,\dots,n$.
For all remaining values of the $j$, set $\omega_j=1$ and $J_j=\emptyset$.

Then, apply Theorem \ref{thm shuffle}. We obtain

\begin{theorem}
\label{thm n-shuffle}
(Shuffle relations) With the above notation, we have
\begin{eqnarray}
\left(\int_{g}^{k',J'}\omega'_{\bf{1}}\dots\omega'_{\bf{k'}}\right)
\left(\int_{g}^{k'',J''}\omega''_{\bf{1}}\dots\omega''_{\bf{k''}}\right)
=\\
=
\sum_{\rho\in Sh(k',k'')}\int_{g}^{k,J}\omega_{\rho({\bf{1}})}\dots\omega_{\rho({\bf{k'+k''}})}.
\end{eqnarray}
\end{theorem}

Let us define
$$
\overline{Sh}(k',k'')=\prod_{i=1}^n \overline{Sh}(k'_i,k''_i),
$$
where the product is over all $i=1,\dots,n$, of shuffles
$\overline{Sh}(k'_i,k''_i)$ as in Definition \ref{def shuffle 1}
Let $j'=(j'_1,\dots,j'_n)$ with $0\leq j'_i\leq k'_1+1$ and
$j''=(j''_1,\dots,j''_n)$ with $0\leq j''_i\leq k'_i+k''_i+1$. For
$\sigma\in\overline{Sh}(k',k'')$, let
$\theta_{\sigma(j')}=\theta'_{j'}$ and
$\theta_{\sigma(j'')}=\theta''_{j''}.$ For all other indecies $j$,
we let $\theta_j=1$. If $\sigma(j')=\sigma(j'')$, (which happens
only at the vertices of the cube, that is, $j'_i=0$ or
$j'_i=k'_i+1$, similarly for $j''$), then we define
$\theta_{\sigma(j')}=\theta'_{j'}\wedge\theta''_{j''}.$

\begin{theorem}
\label{thm n-shuffle 1}
(Shuffle relations) With the above notation, we have
\begin{eqnarray}
\left(\int_{g}^{k',J'}\omega'_{\bf{0}}\dots\omega'_{\bf{k'+1}}\right)
\left(\int_{g}^{k'',J''}\omega''_{\bf{0}}\dots\omega''_{\bf{k''+1}}\right)
=\\
=
\sum_{\rho\in \overline{Sh}(k',k'')}\int_{g}^{k,J}\omega_{\rho({\bf{0}})}\dots\omega_{\rho({\bf{k'+k''+1}})}.
\end{eqnarray}
\end{theorem}

\subsection{Product relations over different membranes}

Let $g'$ and $g''$ be two membranes (images on $n$-cubes into a
manifold $M$) with a common face in direction $i=1$. Let
$\omega'_1,\dots,\omega'_{m'}$ and
$\omega''_1,\dots,\omega''_{m''}$ be differential forms on $M$.
Let $J'$ and $J''$ be multi-indecies as in Subsections 3.1 and 3.2. Let
$\rho'$ and $\rho''$ two $n$-tuples of permutations as in Subsections
3.1 and 3.2.

Consider two iterated integrals
$$I_1=\int_{g'}^{k',J'}\omega'_{\bf{0}}\cdots\omega'_{\bf{m'+1}}
\mbox{ and }
I_2=\int_{g''}^{\rho'',k'',J''}\omega''_{\bf{1}}\cdots\omega''_{\bf{m''+1}}.$$
Let
$$g(t_1,t_2,\dots,t_n)=\left\{
\begin{tabular}{lll}
$g'(2t_1,t_2,\dots,t_n)$ & for & $0<t_1\leq 1/2$\\
\\
$g''(2t_1-1,t_2,\dots,t_n)$ & for & $1/2<t_1<1$
\end{tabular}
\right.
$$

Consider a shuffle product of the integrals $I_1$ and $I_2$ in directions
$i=2,3\dots,n$ and a (non-commutative) product in direction $i=1$.

\begin{definition}
We define $\overline{Sh}^1(k',k'')$ to be a subset of $\overline{Sh}(k',k'')$ such that for
$\rho=(\rho_1,\dots,\rho_n)\in \overline{Sh}^1(k',k'')$ we have that
 $\rho_1$ is the identity. That is $\rho_1$ an automorphism of
$\{0,1,\dots,k'_1+k''_1\}$ sending the first $k'_1+1$ elements to the
first $k'_1$ elements.
\end{definition}
Let also $\omega_j=\omega'_j$ for $j={\bf{1}},\dots,{\bf{m'}}$ and
$\omega_{j+m'}=\omega''_j$ for $j={\bf{1}},\dots,{\bf{m''}}.$

Consider the sum over all shuffles in directions $i=2,3,\dots,n$,
that is, shuffles from the set $\overline{Sh}^1(k',k'')$.  
Using
Theorem \ref{thm shuffle} in directions $i=2,3\dots,n$ and Theorem 2.2 in
direction $i=1$, we obtain the following.
\begin{theorem}
\label{thm n-cutting}
With the above notation we have
\begin{equation}
\label{eq parallel transport}
\int_{g'}^{k',J'}\omega'_{\bf{0}}\cdots\omega'_{\bf{k'+1}}
\int_{g''}^{k'',J''}\omega''_{\bf{0}}\cdots\omega''_{\bf{k''+1}}
=
\sum_{\rho\in\overline{Sh}^1(k',k'')}\int_g^{k,J}\omega_{\bf{0}}\cdots\omega_{\bf{m'+m''+1}}.
\end{equation}
\end{theorem}

\subsection{Connection on higher path space}
Let $P^n=P^1P^{n-1}X=Maps([0,1]^n,X).$ Let $\omega$ be a $p$-form
on $P^{n-1}X$ and let $\theta$ be a $1$-form on $L^{n-1}X$, as
defined in Subsection \ref{subsec forms on higher loops}. We are
going to define a parallel transport of $\omega$ along a path
$\gamma$ in the space $P^{n-1}X$ with respect to a connection
$\nabla=d-\theta$. Note that such a path $\gamma:[0,1]\rightarrow
P^{n-1}X$
 is equivalent to a membrane
$g:[0,1]^n\rightarrow X$.
The parallel transport of $\omega$ along $\gamma$ is
$$\omega
\mapsto
\left(\int_\gamma\omega\right)\theta(\gamma(1))
+
\left(\int_\gamma \omega\cdot\theta\right)\theta(\gamma(1)) +\dots+
\left(\int_\gamma\omega\cdot\theta^{\cdot n}\right)\theta(\gamma(1))+\dots$$

\begin{proposition}
Each of the integrals in the above infinite sum, in the sense of
integrals on differentiable spaces, can be expressed as a finite
linear combination of iterated integrals over an $n$-dimensional
membrane.
\end{proposition}

Let $g_0:[0,1]\rightarrow X$ be an any $(n-1)$-dimensional
membrane, considered as a point on $P^{n-1}X$. Let $\gamma$ be a
path on $P^{n-1}X$. Then $\gamma$ can be considered as a
$n$-dimensional membrane $g:[0,1]^n\rightarrow X$. Let $\omega$
and $\theta$ be forms on $P^{n-1}X$ as defined in Subsection
\ref{subsec forms on higher loops}.

In order to obtain the integral
$\int_\gamma\omega\cdot\theta^{\cdot n}$, we shuffle the
$(n-1)$-dimensional iterated integrals
$\omega,\theta,\dots,\theta$ in each of the directions
$i=1,\dots,n-1$. And we iterate them in direction of $i=n$.

\begin{definition}
We define the set
$\overline{Sh}^n(k',k'',\dots,k'')$
to be the set of shuffles
$\rho=(\rho_1,\dots,\rho_n)\in\overline{Sh}(k',k'',\dots,k'')$.
For $i=1,\dots, n-1,$
$$\rho_i\in Sh(k'_i,k''_i,\dots,k''_i),$$
here we define a shuffle of $l+1$ sets ''$Sh$'' in the same way as a shuffle of $2$ sets.
For $i=n$ let
$\rho_n=id.$
\end{definition}

For $i=1,\dots, n-1$, let
$$k_i=k'_i+k''_i+\dots+k''_i.$$
For $i=n$, let
$$k_n=l+1.$$
Let $dt_n$ be the form for the iteration of the forms $\theta$ and
$\omega$ in direction $\gamma$. Then $dt_n$ occurs exactly once in
the differential form $\omega$. Then $dt_n$ occurs exactly at one
of $\omega'_{\bf{0}},\dots,\omega'_{\bf{k'}}$, say
$\omega'_{j_1}$. Similarly, $dt_n$ occurs exactly once in the form
$\theta$. Then $dt_n$ occurs exactly at one of the forms
$\theta''_{\bf{0}},\dots,\theta''_{\bf{k''}}$, say
$\theta''_{j_2}$.

Then, for $j=1,\dots,m'$,
$J_j=J'_j$ for $j\neq j_1,$
and
$J_{j_1}=J'_{j_1}\cup\{n\};$
for $j=m'+(a-1)m''+1,\dots,m'+am''$, $a=1,\dots,l-1$, we define
$J_j=J'_{j-m'-am''}$ for $j-m'-am''\neq j_2$
and
$J_{j_2+m'+am''}=J_{j_2}\cup\{n\};$
for $j=m'+(l-1)m''+1,\dots,m'+lm''$, we define
$J_j=J'_{j-m'-lm''}.$

\begin{theorem}
\label{thm higher connection} With the above notation, a parallel
transport of a form $\omega$ on $P^{n-1}X$ with respect to a
connection $\theta$ on $P^{n-1}X$, given at a finite step by
$$\omega\mapsto\left(\int\omega\cdot\theta^{\cdot n}\right)\theta$$
is given by a finite linear combination of iterated integrals on a
membrane over the manifold $X$. Explicitly, it is achieved in the
following way. Let
$$\omega_{g_0}=\int_{g_0}^{k'_0,J'_0}\omega'_{\bf{0}}\cdots\omega'_{\bf{k'+1}}$$
and
$$\theta_{g_0}=\int_{g_0}^{k''_0,J''_0}\theta''_{\bf{0}}\cdots\theta''_{\bf{k''+1}}.$$
Then
$$\left(\int_\gamma \omega\cdot\theta^{\cdot(l-1)}\right)\theta
=
\sum_{\rho\in\overline{Sh}^n(k',k'',\dots,k'')}
\int_g^{k,J}
\omega_{\bf{0}}\dots\omega_{\bf{k+1}}.$$
\end{theorem}
\subsection{De Rham structure on higher homotopy groups}

Recall a notation from Chen. Let $\hat{C}^*(X)$ be a chain complex
of those smooth singular simplices of $X$ that map the
$1$-skeleton of the simplex to the base point $x_1$ of $X$. Let
$L_{x_0}X$ denotes the space of smooth maps from $(S^1,*)$ to $(X,x_1)$.
Let $UX$ denote the universal cover of $X$. Consider the space
$UL_{x_{n-1}}U\dots L_{x_1}UX$.

Using the iterated integrals on membranes from the definitions in
Subsections 3.1 and 3.2 we can express a connection on the
$(n-1)$-loop space by such integrals (Theorem \ref{thm higher
connection}). Then Chen's theorem about the fundamental group,
Theorem \ref{thm pi1}, leads to the following.

\begin{theorem}
Let $X$ be a manifold. Let $W$ be a
vector space of closed $1$-forms on $L_{x_{n-1}}U\dots L_{x_1}UX$, representing
cohomology classes, with basis $\{w_1\dots,w_m\}$. Let
$V=\{v_1,\dots,v_l\}$ be a basis for $W\wedge W\subset
\Lambda^2(X)$. We define the coefficients $c_{ijk}$ so that
$$w_i\wedge w_j=\sum c_{ijk}v_k.$$
Then there exists an epimorphism of Lie algebras
$$gr \pi_n(X)\otimes_\Z\R\rightarrow L(X)/(N),$$
where $(N)$ is the ideal generated by $\sum_{i,j}c_{ijk}[X_i,X_j]$ for $k=1,\dots,l$.
\end{theorem}
\proof
Clearly $\pi_1(L_{x_{n-1}}U\dots L_{x_1}UX,x_n)=\pi_n(X,x_1)$ and $L_{x_{n-1}}U\dots L_{x_1}UX$  is a differentiable space,
when $X$ is a manifold. Then this Theorem is a consequence of Theorem 2.11.

It will be interesting to examine the kernel of the above map. For that pupose one can examine the kernel on the level of the fundamental group of a differetiable space and apply the results of Chen from \cite{Ch2}.

\subsection{De Rham complex on higher loop spaces}

Using the iterated integrals on membranes from the definitions in
Subsections 3.1 and 3.2 we can express a connection on the
$(n-1)$-loop space by such integrals (Theorem \ref{thm higher
connection}).

Denote by $LU\Lambda^*(X)$ the de Rham complex on $L_{x_1}UX$ defined via
iterated integrals, when $\Lambda^*(X)$ is the de Rham complex
on $X$. By induction on the dimension, we define $(LU)^n\Lambda$
the cochain complex of differential form, in terms of iterated
integrals on $n$-dimensional membranes, which are maps from a
torus such that $n-1$ of the generators of the homology of the
torus map to homotopy trivial paths.

Then Chen's theorem about the fundamental group,
Theorem \ref{thm 1loop}, leads to the following.
\begin{theorem}
Assume that

(a) $H_*(U(LU)^{n-1}X,\R)$ is of finite type;

(b) If $X$ is a manifold then the chain complex
$\hat{C}_*(UL_{x_{n-1}}U\dots L_{x_1}UX)$ is quasi-isomorphic to
$C_*(UL_{x_{n-1}}U\dots L_{x_1}UX)$.

(c) $H^*(U(LU)^{n-1}\Lambda(X))\cong H^*(U(LU)^{n-1}X,\R)$;

Then $H^*((LU)^n\Lambda^*(X))=H^*((L_{x_n}U\dots L_{x_1}UX)$.
\end{theorem}

For the Definition of a chain $\hat(C)_*$ see Subsection 2.6.
\begin{conjecture}
Conditions (a), (b) and (c) are satisfied for smooth complex quasi-projective algebraic varieties.
\end{conjecture}

\section{Final remarks}
In Subsections 3.6 and 3.7, we considered the space
$L_{x_{n-1}}U\dots L_{x_1}UX$ we a higher loop space. It has the
homotopy type of an $(n-1)$fold loop space. However the space is
$L_{x_{n-1}}U\dots L_{x_1}UX$ is a map of $(n-1)$-dimensional
torus $T^{n-1}$ to $X$, so that $n-2$ of the generators of
$H_1(T^{n-1})$ of the torus map to  trivial loops on $X$, up to
homotopy. Let us define inductively the space of maps from an $n$-torus to the manifold $X$
$$Maps((T^{n+1},\dot{T}^{n+1}),(X,x_1,\dots,x_{n+1}))=Maps((S^1,*),Maps((T^{n},\dot{T}^{n}),(X,x_1,\dots,x_{n})),$$
where  
$$x_{n+1}\in Maps((T^n,\dot{T}^n),(X,x_1,\dots,x_n))$$ and
$$Maps((T^{1},\dot{T}^{1}),(X,x_1))=Maps((S^1,*),(X,x_1)).$$
Let
$$Maps((T^\infty,\dot{T}^\infty),(X,\{x_n\}))$$
be the space of maps  from an infinite dimensional torus, which factor through a map from a finite dimensional torus as above.
\begin{conjecture} (a) The space
$Maps((T^\infty,\dot{T}^\infty),(X,\{x_n\}))$, together with all
its quotients $Maps((T^n,\dot{T}^n),(X,x_1,\dots,x_n))$ capture
the homotopy type of $X$ (together with the Postnikov
tower).

(b) Iterated integrals over membranes on a smooth complex quasi-projective algebraic variety $X$
capture the rational homotopy type of $X$ together with the Postnikov
tower up to tensoring with the rational numbers.
\end{conjecture}
{\bf{Reason for the Conjecture:}}
(a) Form Subsections 3.6, it is clear that
$$Maps((T^n,\dot{T}^n),(X,x_1,\dots,x_n))$$
capture $\pi_i(X)$, for $i=1,\dots,n$, under mild
assumptions. Assume that
$$Maps((T^n,\dot{T}^n),(X,x_1,\dots,x_n)),$$ captures also
rational Postnikov tower up to the $(n+1)$-dimensional class.
Consider the map from an $(n+1)$-dimensional torus,
$Maps((T^{n+1},\dot{T}^{n+1}),(X,x_1,\dots,x_{n+1}))$, as the torus is
parametrized by a $(n+1)$-dimensional cube. Then the $(n+2)$-dimensional cohomology
class in the Postnikov tower can be realized as the boundary of a
variation of a function in
$$Maps((T^{n+1},\dot{T}^{n+1}),(X,x_1,\dots,x_{n+1})).$$
Two of the boundary components are maps from $T^{n+1}$ and the
remaining components are variations of maps from
$$Maps((T^n,\dot{T}^n),(X,x_1,\dots,x_n)).$$
Note that the last space caries information about $\pi_1(X),\dots,\pi_n(X)$ together with the Postnilov tower up to the $(n+1)$-dimensional cohomology class. 
This variation gives the $(n+2)$-dimensional cohomology class. More importantly, it gives the induction step for constructing the Postnikov tower.

(b) From Subsection 3.6, it is clear that iterated integrals on membranes capture rational homotopy of a smooth complex quasi-projective algebraic variety $X$. On the de Rham side 'variation' means the differential $d$ on functions on higher loop spaces, realized as iterated integrals on membranes. 
Not surprisingly, the above described relation is a type of Stokes theorem. 

{\bf{Acknowledgements:}} I am indebted to Manin for his interest and encouragement
regarding my construction of iterated integrals on membranes
\cite{M2} and for his inspiring presentation of \cite{M}. 
In \cite{M2} he pointed out that such integrals could
be used to describing higher loop spaces. I would like to thank R.
Brown for the inspiring talk that he gave at Durham University,
based on which I constructed the first version of iterated
integrals on membranes \cite{BH}. I would like to thank Goodwillie
for pointing out some errors in the first preprint on iterated
integrals on membranes \cite{ModSymb} and
Sinha for his interest in this project, which started with
\cite{DH}, and for the discussions we had.

Finally, I gratefully acknowledge the  financial
support and good working environment provided by Max Planck Institute for Mathematics, Durham University, Brandeis University, Tubingen University and Washington University in St Louis during the period of working on this paper.

\renewcommand{\em}{\textrm}

\begin{small}

\renewcommand{\refname}{ {\flushleft\normalsize\bf{References}} }
    
\end{small}
Ivan Horozov\\
\begin{small}
Department of Mathematics\\
Washington University in St Louis\\
One Brookings Drive
St Louis, MO 63130\\
USA\\
horozov@math.wustl.edu

\end{small}
\end{document}